\documentclass[12pt,a4paper]{article}

\usepackage[utf8]{inputenc}
\usepackage[english]{babel}
\usepackage{amsmath}
\usepackage{amsfonts} % to enable \mathbb
\usepackage{amssymb}
\usepackage{makeidx}
\usepackage{graphicx}
\usepackage{lmodern}
\usepackage{braket} % to enable \Set
\usepackage{caption} % to enable subfigure
\usepackage{subcaption} % to enable subfigure
\usepackage{algorithm}
\usepackage[hyphens]{url} % 
\usepackage[top=25mm,bottom=25mm]{geometry}

\newtheorem{theorem}{Theorem}
\newtheorem{remark}[theorem]{Remark}

\title{How much testing and social distancing is required to control COVID-19? 
    Some insight based on an age-differentiated compartmental model%
	\thanks{This work was funded by the Federal Ministry of Education and Research (BMBF; grants 05M18EVA and 05M18SIA).}} %

\author{%
	Sara Grundel%
	\thanks{Max Planck Institute for Dynamics of Complex Technical Systems, Magdeburg, Germany
		(\texttt{grundel@mpi-magdeburg.mpg.de}, \texttt{ritschel@mpi-magdeburg.mpg.de}).
	}
	\and
	Stefan Heyder%
	\thanks{Technische Universit{\"{a}}t Ilmenau, Ilmenau, Germany, Institute for Mathematics 
		(\texttt{stefan.heyder@tu-ilmenau.de}, \texttt{thomas.hotz@tu-ilmenau.de}, \texttt{philipp.sauerteig@tu-ilmenau.de}, \texttt{karl.worthmann@tu-ilmenau.de}).
		}
	\and
	Thomas Hotz\footnotemark[3]
	\and
	Tobias K. S. Ritschel\footnotemark[2]
	\and
	Philipp Sauerteig\footnotemark[3]
	\and
	Karl Worthmann\footnotemark[3]
}

\begin{document}
	% Title, etc.
	\maketitle
	
	% Abstract, 
	\begin{abstract}
In this paper, we provide insights on how much testing and social distancing is required to control COVID-19. To this end, we develop a compartmental model that accounts for key aspects of the disease: 1)~incubation time, 2)~age-dependent symptom severity, and 3)~testing and hospitalization delays; the model's parameters are chosen based on medical evidence, and, for concreteness, adapted to the German situation. Then, optimal mass-testing and age-dependent social-distancing policies are determined by solving optimal control problems both in open loop and within a model predictive control framework. We aim to minimize testing and/or social distancing until herd immunity sets in under a constraint on the number of available intensive care units. We find that an early and short lockdown is inevitable but can be slowly relaxed over the following months.
%}
\end{abstract}
	
	% Main content
	\section{Introduction}
The \emph{severe acute respiratory syndrome coronavirus 2} (SARS-CoV-2) is a strain of coronavirus which causes the respiratory illness \emph{coronavirus disease 2019} (COVID-19).
In 2020, on March 11\textsuperscript{th}, the World Health Organization (WHO) declared the outbreak of SARS-CoV-2 a pandemic~\cite{WHO20}.
Due to the novelty of the virus, there was (and, at the time of submitting this manuscript, still is) significant uncertainty concerning the severity and mortality of COVID-19.
Furthermore, as of October 2020, no vaccine has completed the trials necessary for approving widespread use~\cite{LeCRetal20}.
Therefore, many countries are enforcing nonpharmaceutical countermeasures \cite{IMF20c, IMF20b}, e.g. 1)~social distancing, 2)~increased public hygiene, 3)~travel restrictions, 4)~self-isolation (quarantine), and 5)~population-wide mass testing for SARS-CoV-2 infection.
However, enforcing these countermeasures for long periods of time can have severe economic and social consequences, both at the national and the global scale~\cite{OECD20}. Therefore, there is a need for identifying economic strategies for simultaneously relaxing the countermeasures and containing the pandemic.

Model-based decision support systems can be used for exactly this purpose. They use \emph{predictive} models to assist decision makers in identifying and evaluating candidate strategies (e.g. \cite{AroGM20}).
In particular, given a dynamical model of the spread of SARS-CoV-2, economically optimal (open-loop) mitigation strategies can be identified by solving optimal control problems (OCPs) over several months or even years. A key advantage of this approach is that it can directly account for constraints, e.g. related to the capacity of public healthcare systems. However, given the uncertainty surrounding SARS-CoV-2 and COVID-19, it is advisable to implement the optimal mitigation strategies in closed-loop, i.e. to repeatedly update the strategies when new data becomes available. This is referred to as model predictive control (MPC)~\cite{RawMD19} and is an established method for advanced process control~\cite{QinB03}. The predictive capabilities of the underlying model are crucial for the efficacy of the resulting mitigation strategy, and a common challenge is to identify suitable model parameters.

Epidemics are often modelled using deterministic \emph{compartmental models}~\cite{Het00}, e.g. the classical SIR model, where individuals are either susceptible, infectious, or removed, or the SEIR model which, additionally, takes the incubation time into account.
Optimal control of compartmental models was already an active research topic before the SARS-CoV-2 pandemic (see \cite{ShaM17} for a review). In particular, optimal control of SIR models has been considered, e.g. for arbitrary social interaction models~\cite{Beh00} and to identify time-optimal mitigation strategies~\cite{BolBS17, HanD11}.
Optimal control of more complex models has also been considered. For instance, Fischer et al.~\cite{FiscChud19} consider optimal control of a model with two species, Bussell et al.~\cite{BusDC19} demonstrate the importance of closed-loop mitigation strategies (i.e. of incorporating feedback), and Watkins et al.~\cite{WatNP20} consider MPC of stochastic compartmental models. %
In~\cite{EstALetal20} the authors determine control strategies to maintain hard infection caps in a disease-vector model based on the theory of barriers. This approach, however, exploits the low dimensionality of the model. Application of these techniques to complex compartmental models, therefore, requires model order reduction. 
%}

In response to the SARS-CoV-2 pandemic, many researchers have presented optimal control strategies, for instance
based on Pontryagin's maximum principle (e.g. \cite{KanK20, PerE20, ZamSNetal20}). These strategies typically involve 1)~extended SIR or SEIR models, 2)~nonpharmaceutical countermeasures (often social distancing), and 3)~minimization of the number of infected as well as the economic cost of the countermeasures (and often other quantities as well, e.g. the number of deaths). Furthermore, they rarely satisfy hard constraints, for instance related to health care or testing capacities.
In the following, we highlight some of the key developments in decision support for SARS-CoV-2 mitigation based on optimal control.
Gondim and Machado~\cite{GonM20} use a model with three age groups to compute optimal quarantine strategies (for susceptible individuals) which minimize the number of infected and the cost of quarantining. Bonnans and Gianatti~\cite{BonG20} compute social distancing strategies based on a model with a continuous age structure. Here, the strategies minimize a combination of 1)~the number of deaths, 2)~the peak number of hospitalized, and 3)~the cost of social distancing. Similarly, Richard et al.~\cite{RicACetal20} present optimal social distancing strategies based on a model with a continuous age and infection duration structure, which minimize the number of deaths and the cost of social distancing.
Morato et al.~\cite{MorBCetal20} compute on-off (also called bang-bang) social distancing strategies which minimize 1)~the number of symptomatic infectious people and 2)~the duration of the social distancing policies, subject to constraints on intensive care unit (ICU) occupancy. They use extended SIR models.
Carli et al.~\cite{CarCEetal20} use MPC to compute social distancing and travel restriction strategies for an extended multi-region SIR model, minimizing the cost of the countermeasures and preventing an overload on the hospitals.
K{\"{o}}hler et al.~\cite{KohSjetal20} use MPC to minimize the number of fatalities caused by COVID-19, subject to constraints on the economic cost of social distancing. They take a modified SIDARTHE model~\cite{GioBBetal20} as basis and use interval arithmetic in the MPC to propagate model uncertainties.
Finally, Tsay et al.~\cite{TsaLSetal20} use MPC to minimize the cost of social distancing and testing, subject to an upper bound on the peak number of infectious people who have been tested positive. They use the unscented Kalman filter to estimate the noisy state variables of an extended SEIR model.

In this work, we address some of the key questions that decision makers involved in the mitigation of the SARS-CoV-2 pandemic are facing: 1)~Is mass testing alone sufficient to avoid overloading of ICUs? 2)~If not, how much social distancing is then required? 3)~How much can social distancing measures be reduced by targeting specific age groups? 4)~How do strategies obtained by short and long-term planning differ? 5)~What are the benefits of increasing the daily testing capacity or the ICU capacity? Here, the limited ICU capacity is considered as an example for constraints imposed by the health care system or political considerations. Of course, different constraints such as limited personnel for contact tracing could be incorporated as well.

We address the above questions by proposing a novel compartmental model and using optimal control as well as MPC to compute open and closed-loop social distancing and testing strategies. The model contains three age groups, and it accounts for several of the key challenging characteristics of COVID-19, i.e. 1)~the incubation time, 2)~different levels of symptom severity depending on age, 3)~delay of testing results (and the following self-isolation), and 4)~delay of hospitalization. Furthermore, we choose values of the epidemiological model parameters based on the current state of knowledge in order to ensure that our numerical results match reality. For concreteness, we use the COVID-19 outbreak in Germany to determine parameters depending on demographics and the health care system. However, we expect our conclusions to carry over to outbreaks in other developed countries as well.

The remainder of this paper is structured as follows. In Section~\ref{sec:model}, we describe the novel compartmental model of the SARS-CoV-2 outbreak in Germany, and in Section~\ref{sec:parameters}, we motivate our choice of model parameters. In Section~\ref{sec:optimal_control},
we demonstrate that optimal control can be used as a decision support tool based on the proposed model, 
and we conclude the paper in Section~\ref{sec:conclusions}.

\newcommand{\rS}{\mathcal S}
\newcommand{\rI}{\mathcal I}
\newcommand{\rR}{\mathcal R}
\newcommand{\ind}{\mathbf 1}

\section{Modelling  Pandemics}\label{sec:model}
In this section, we propose a dynamical model tailored to \mbox{COVID-19}. The aim is to be able to evaluate the effect of population-wide mass testing (in combination with quarantine) and social distancing measures on the development of the pandemic. To this end, we extend the well-known SIR model.

\subsection{Interpretation of deterministic compartmental models}
We start with an illustration of the connection between 1)~infectious disease models based on randomly acting individual agents and 2)~their approximation by ordinary differential equation compartmental models. 
This exposition will highlight the interpretation and conversion of parameters when moving from a random to a deterministic model. 
For simplicity, we consider the classical SIR model in this subsection. 
However, the connection, especially the interpretation of parameters, is similar for more complex models such as the one described in Section~\ref{sec:novel:seitphr:model}.

Consider a population of $ n_{\mathrm{pop}} $ individuals or agents each being either \textit{susceptible}, \textit{infectious} or \textit{removed}.
At time $ t \in [0, \infty) $ denote the (random) set of susceptibles by $ \rS_t $, the set of infectious by $ \rI_t $ and the set of removed by $ \rR_t $.
Time is modeled continuously and measured in days.

We assume a homogeneous population with contacts between agents $ a $ and $ b $ following a Poisson process with intensity $ \lambda $ which does not depend on the agents considered.
Infections occur randomly upon contact with a fixed probability $ \alpha $ if one of the agents is susceptible and the other infectious.
Thus, potentially infectious contacts also follow a Poisson process with respective intensity $ \alpha\lambda $.
Similarly, we model other events, in this simple model only recoveries, to occur according to a Poisson process. 
This implies that the time an agent spends in the infectious compartment is exponentially distributed with rate $ \eta $, say, which we also assume to be the same for each agent (see~\cite{OguP20} for models where these quantities follow other distributions).

We denote by $ S(t) = \mathbf E \frac{\lvert \rS_t \rvert}{n_\mathrm{pop}} , I(t) = \mathbf E  \frac{\lvert \rI_t \rvert}{n_\mathrm{pop}} $ and $ R(t) = \mathbf E \frac{\lvert \rR_t \rvert}{n_\mathrm{pop}}  $ the expected share of the population which are susceptible, infectious and removed, respectively.
Since for large $ n_{\mathrm{pop}} $ the change of $\frac{\lvert \rS_t \rvert}{n_\mathrm{pop}}$ over a short time interval can, due to the law of large numbers, be well approximated by its expectation, $ S(t) $ will provide a sufficient approximation of $ \frac{\lvert \rS_t \rvert}{n_\mathrm{pop}} $ over the finite time horizon considered for a country the size of Germany.
By the same argument, $ I(t) $ and $ R(t) $ approximate $ \frac{\lvert \rI_t \rvert}{n_\mathrm{pop}} $  and $\frac{\lvert \rR_t \rvert}{n_\mathrm{pop}}$, respectively, sufficiently well.

If $ a $ is susceptible he will transition to the infectious compartment upon having an infectious contact. 
At a fixed time $ t $ with $ a \in \rS_t $, there are two possible sources of infection for $ a $: either some $ b \in \rI_t $ which is already infectious or some $ c \in \rS_t $ which will become infectious himself at some later time.

To determine the probability that $ b $ infects $ a $ in the time frame $ (t, t + \Delta t] $, we analyze two competing events:
The first is an infectious contact between $ a $ and $ b $, and the second is $ b $'s recovery from the infectious state.
Both events happen independently of one another with exponentially distributed time of occurrence, the first with rate $ \alpha \lambda $ and the second with rate $ \eta $. 
Thus the first time of occurrence of one of these is again exponentially distributed with rate $ \alpha\lambda + \eta $ and the probability that the first occurrence is an infectious contact is $ \tfrac {\alpha\lambda}{\alpha\lambda + \eta} $. 
In total
\begin{align*}
    & \mathbf P \big( b \text{ infects } a \text{ in } (t, t + \Delta t] ~\big\vert~ a \in \rS_t, b \in \rI_t \big) \\
    = & \big( 1 - \exp(-(\alpha\lambda + \eta)\Delta t) \big) \frac{\alpha\lambda}{\alpha\lambda + \eta} \\ 
    = & \alpha\lambda\Delta t + o(\Delta t).
\end{align*}

For $ c \in \rS_t$ to infect $ a $ in $ (t, t + \Delta t] $, $ c $ has to become infectious himself before he in turn can infect $ a $. 
This happens only with probability $ o((\Delta t)^2) $ and can, thus, be neglected in the following calculations.
In total $ a $ is moved out of the susceptible compartment with probability
\begin{align*}
    \mathbf P \Big(a \notin \rS_{t + \Delta t} ~\big|~ \rS_t, \rI_t \Big) &= 
    1 - \prod_{b \in \rI_t} \Big( 1 - \alpha \lambda \Delta t  - o(\Delta t)\Big) \prod_{c \in \rS_t} \Big(1 - o\big((\Delta t)^2\big)\Big)\\
                                                                              &= -\alpha\lambda \Delta t \lvert \rI_t\rvert +  o(\Delta t).
\end{align*}
Approximating $ \frac {\lvert\rS_t\rvert}{n_\mathrm{pop}} $ and $ \frac{\lvert \rI_t\rvert}{n_\mathrm{pop}} $  by $ S(t) $ and $ I(t) $ using the law of total expectation yields 
\begin{align*}
    \frac{S(t + \Delta t) -  S(t)}{\Delta t} &\approx \frac{1} {n_\mathrm{pop} \Delta t}\sum_{b \in \rS_t} -\alpha \lambda \Delta t \lvert \rI_t\rvert + o(\Delta t) = -\frac{\alpha\lambda}{n_\mathrm{pop}} \lvert \rS_t \rvert \lvert \rI_t\rvert + o(1)\\ 
                                             &\approx - \alpha \lambda n_\mathrm{pop} S(t) I(t) + o(1).
\end{align*}
As we assume the time from infection to removal to be exponentially distributed with rate $ \eta $, a similar but more straight-forward calculation reveals
\begin{align*}
    \frac{R(t + \Delta t) - R(t)}{\Delta t} = \eta I(t) + o(1),
\end{align*}
where $ \eta^{-1} $ is the mean stay of a single agent in the infectious compartment.
We now set $ \beta = n_\mathrm{pop} \alpha \lambda $, which can be interpreted in this model as the daily amount of (potentially) infectious contacts a single agent has.

Since $ S(t) + I(t) + R(t) = 1 $ for all $ t $, we obtain the following system of ODEs:
\begin{equation}
    \begin{aligned}\label{eq:sir}
        \dot S(t) &= - \beta S(t) I(t) \\
        \dot I(t) &= \beta S(t)I(t) - \eta I(t)\\
        \dot R(t) &= \eta I(t).
    \end{aligned}
\end{equation}
To determine suitable parameter values for $ \beta $ and $ \eta $ in this simplistic model, we reiterate that these are best thought of in the probabilistic setting.
For the coefficients of the linear terms on the right-hand side, the interpretation is straight-forward: it is the rate of the exponential distribution underlying the time until an agent leaves the compartment.
Its inverse is the mean stay in this compartment. 

For coefficients of interaction (product) terms, here $ \beta $, the interpretation is the rate at which an agent in the first compartment causes other agents to leave the second compartment. 
In our setting, this is the daily amount of infections one infectious agent causes which can readily be seen from the definition of $ \beta $.
See Section \ref{sec:parameters} for a more detailed discussion of the parameter values we use in our model.
As above mentioned, these interpretations for the parameters carry over in a straight-forward manner to more sophisticated models such as the one considered in the following.

\subsection{A Compartmental Model for COVID-19}\label{sec:novel:seitphr:model}
\begin{figure}[htbp!]
	\centering	\includegraphics[width=0.9\textwidth]{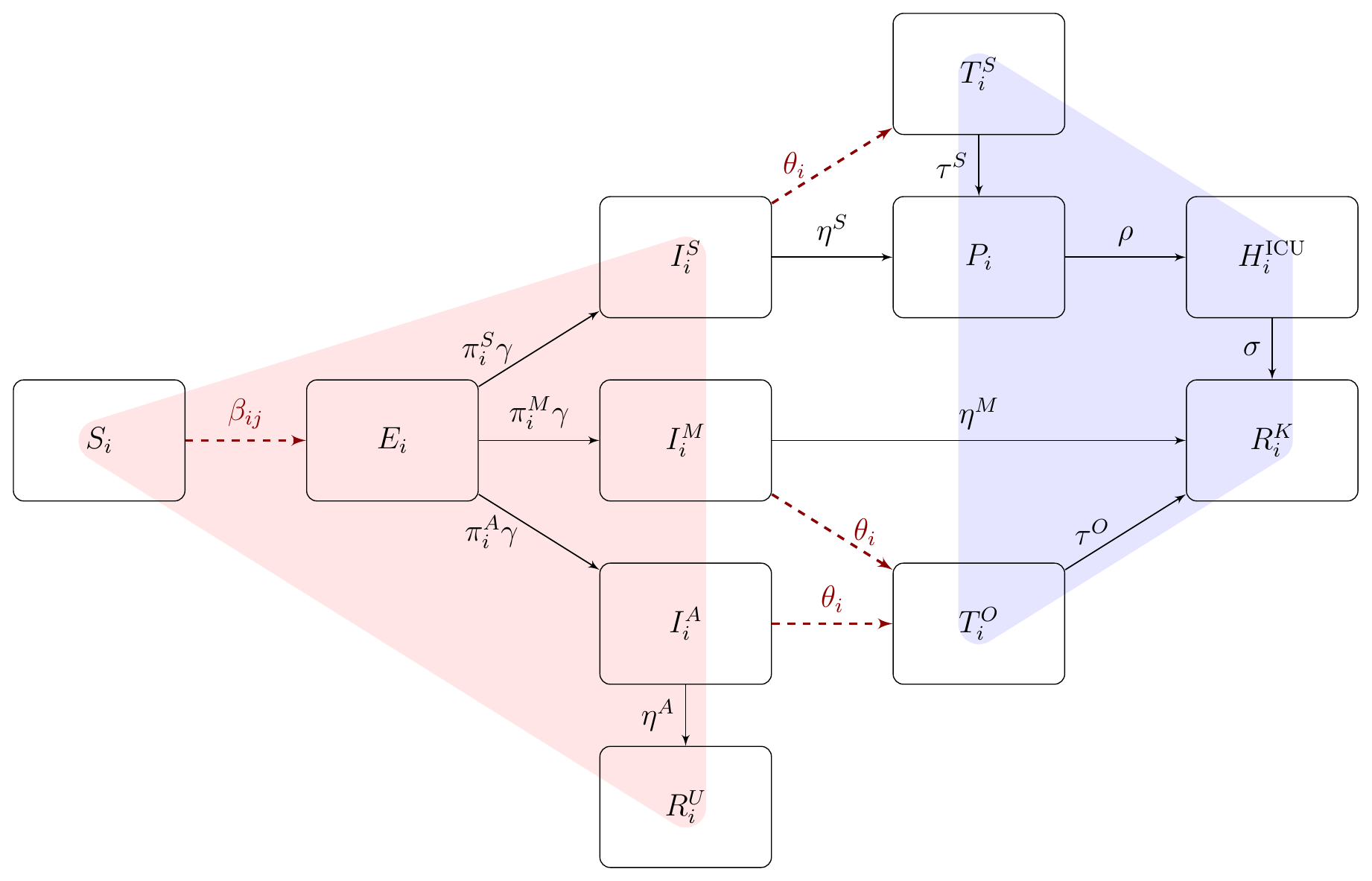}
	\caption{Flow of the SEITPHR model for one age group. The controls are indicated with dashed red edges. Unreported compartments are highlighted by the left red triangle, while tested and detected compartments are highlighted by the right blue trapezium.}
	\label{fig:flow_model}
\end{figure}

The SIR model provides a good starting point to study the dynamics of pandemics. However, due to its simple structure it is not suited to model the \mbox{COVID-19} pandemic adequately. In particular it does not include hospitalization, age-specific disease progressions and interventions. Therefore, we extend the SIR model in three ways. 
\begin{enumerate}
	\item We introduce eight additional compartments. In detail, we take into account that people can be infected, but not yet infectious. We call them \emph{exposed} (or latent) and denote the compartment by~$E$, see also~\cite{Het00}. Moreover, we split the infectious compartment into three depending on how the course of the infection will be. We distinguish between \emph{severe} cases~$I^S$ that are going to need intensive care, i.e. they will move to $H^\mathrm{ICU}$ at some point in time; \emph{mild} cases~$I^M$ that are going to visit a physician and hence be quarantined, i.e. removed; and \emph{asymptomatic} cases~$I^A$ that might recover without being detected. Furthermore, we incorporate the possibility of being \emph{tested} but not yet detected by introducing the compartments~$T^S$ (severe) and~$T^O$ (other). We assume that the patients with severe cases will visit a physician at some point before being sent to an ICU. To this end, we introduce~$P$ as a pre-ICU compartment which comprises isolated patients at home or on a regular hospital ward. Moreover, we split the compartment of removed people into \emph{known} and \emph{unknown} cases $R = R^U + R^K$. %
	\item Each compartment is further divided into $n_\mathrm{g}$ groups, $n_\mathrm{g} \in \mathbb{N}$, depending on the age of a subject in order to study how these groups affect each other. %
	\item Social distancing and hygiene measures affect the contact rate as well as the transmission probability. Therefore, $\beta$ can be used as a time-dependent control input $ \beta(t) $. 
\end{enumerate}
The resulting SEITPHR model reads as
\begin{subequations}\label{eq:model}
	\begin{align}
	\dot S_i(t) & \; = \; -\sum_{j=1}^{n_\mathrm{g}} \beta_{ij}(t) S_i(t) \left[ I_j^S(t) + I_j^M(t) + I_j^A(t) + T_j^S(t) + T_j^O(t) \right] \label{eq:model_S} \\
	\dot E_i(t) & \; = \; \sum_{j=1}^{n_\mathrm{g}} \beta_{ij}(t) S_i(t) \left[ I_j^S(t) + I_j^M(t) + I_j^A(t) + T_j^S(t) + T_j^O(t) \right] - \gamma E_i(t) \\
	\dot I_i^S(t) & \; = \; \pi_i^S \gamma E_i(t) - (\eta^S + \theta_i(t)) I_i^S(t) \\
	\dot I_i^M(t) & \; = \; \pi_i^M \gamma E_i(t) - (\eta^M + \theta_i(t)) I_i^M(t) \\
	\dot I_i^A(t) & \; = \; \pi_i^A \gamma E_i(t) - (\eta^A + \theta_i(t)) I_i^A(t) \\
	\dot T_i^S(t) & \; = \; \theta_i(t) I_i^S(t) - \tau^S T_i^S(t) \\
	\dot T_i^O(t) & \; = \; \theta_i(t) \left[ I_i^M(t) + I_i^A(t) \right] - \tau^O T_i^O(t) \\
	\dot P_i(t) & \; = \; \eta^S I_i^S(t) + \tau^S T_i^S(t) - \rho P_i(t) \\
	\dot H_i^\mathrm{ICU}(t) & \; = \; \rho P_i(t) - \sigma H_i^\mathrm{ICU}(t) \\ 
	\dot R_i^K(t) & \; = \; \eta^M I_i^M(t) + \tau^O T_i^O(t) + \sigma H_i^\mathrm{ICU}(t) \\
	\dot R_i^U(t) & \; = \; \eta^A I_i^A(t),
	\end{align}
\end{subequations}
where the subscript $i \in \{1,2,\ldots,n_\mathrm{g}\}$ denotes the age group in ascending order. We enforce $\sum_{i=1}^{n_\mathrm{g}} N_i = 1$, where $N_i$ denotes the relative size of age group~$i$. We assume a mean incubation time~$\gamma^{-1}$ independent of both the course of infection and the age of the patient. However, depending on the age, there are different probabilities $\pi_i^S$, $\pi_i^M$, and $\pi_i^A$ for the three courses of infection, where $\pi_i^S + \pi_i^M + \pi_i^A = 1$ for all $i$. Similar to the SIR model~\eqref{eq:sir}, the parameters denoted by~$\eta$ correspond to people being removed from the system, i.e. $\eta^S$ and $\eta^M$ denote those who visit a physician and, therefore, are put into quarantine immediately, while $\eta^A$ represents unreported recovery. We denote the total number of susceptibles and unreported cases by $U_i = S_i + E_i + I_i^S + I_i^M + I_i^A + R_i^U$. The control input $\theta_i : \mathbb{R}_{\geq 0} \to \mathbb{R}_{\geq 0}$ describes the rate of those being tested per day, %. %
where tests are distributed uniformly at random among all individuals in~$U_i$. %
In addition, symptomatic cases who visit physicians are assumed to be tested as well. Therefore, the total number of tests at time $t \geq 0$ is given by
\begin{align}
	T^\mathrm{tot}(t) \; = \; n_\mathrm{pop} \cdot \left( \sum_{i=1}^{n_\mathrm{g}} \theta_i(t) U_i(t) + \eta^S I_i^S(t) + \eta^M I_i^M(t) \right). \notag
\end{align}
Note that testing does not affect the state of non-infectious subjects. Parameters~$\tau^S$ and~$\tau^O$ denote the rate from being tested positive to being detected, and hence being put into quarantine. Furthermore, $\rho$ is the rate from pre-hospital quarantine to hospitalization and $\sigma$ from hospitalization to being reportedly removed, i.e. $\sigma$ incorporates both mortality and recovery rate of hospitalized patients. The basic structure of the SEITPHR model~\eqref{eq:model} is depicted in Figure~\ref{fig:flow_model}.

For a concise notation we stack the state vectors into $x = (x_1, \ldots x_{n_\mathrm{g}})$
and the controls into $u = (\beta, \theta)$, where 
\begin{align*}
	x_i = (S_i, E_i, I_i^S, I_i^M, I_i^A, T_i^S, T_i^O, P_i, H_i^\mathrm{ICU}, R_i^U, R_i^K),
\end{align*} 
and $\beta = (\beta_{ij})_{i,j=1}^{n_\mathrm{g}}$ with $\beta_{ij} : \mathbb{R}_{\geq 0} \to \mathbb{R}_{\geq 0}$ %
and $\theta = (\theta_1, \ldots, \theta_{n_\mathrm{g}})$. Similarly, we denote $\pi \in \mathbb{R}^{3 n_\mathrm{g}}$, $\tau \in \mathbb{R}^2$, and $\eta \in \mathbb{R}^3$. Thus, we write system~\eqref{eq:model} as 
\begin{align}
	\dot{x}(t) \; = \; f(x(t),u(t),p), \label{eq:xdot=f(x,u,p)}
\end{align}
where $p = (\pi, \eta, \tau, \rho, \sigma, \gamma) \in \mathbb{R}^{3n_\mathrm{g} + 8}$ collects all parameters. Furthermore, we introduce the initial condition $x(0) \; = \; x^0$ for some $x^0 \in \Omega$, where 
\begin{align}
	\Omega \; = \; \Set{x \in \mathbb{R}_{\geq 0}^{11 \cdot n_\mathrm{g}} | \sum_{j=1}^{11 \cdot n_\mathrm{g}} x_j = 1} \notag
\end{align} 
denotes the set of possible states. Note that~$\Omega$ is forward invariant under~\eqref{eq:xdot=f(x,u,p)}, i.e. if $x^0 \in \Omega$ then $x(t) \in \Omega$ for all $t \geq 0$.

\section{Parameters}\label{sec:parameters}
Before we present our choices for the parameters of model \eqref{eq:model}, let us reiterate that some of the parameters of our model depend on age.
We indicate this dependence by an appropriate index which we drop if the parameter is constant across age-groups.
For example, $ \pi_i^S $ is the age-dependent probability of having a severe course of disease while we assume $ \eta^A $, the rate with which asymptomatic cases recover, to be age-independent.

\paragraph{\boldmath $N_i$}
We use data on the population size of Germany at the end of $ 2019 $ from the GENESIS-Online Database of the DESTATIS \cite{Genesis2019}. 
The first age group consists of individuals aged younger than $ 15 $ years, the second of those older than $ 15 $ but younger than $ 60$ years, and the last comprises all individuals older than $ 60 $ years.
These groupings result in proportions $ N_1 = 0.14, N_2 = 0.58 $ and $ N_3 = 0.28 $.

\paragraph{\boldmath $\beta^{\mathrm{0}}_{ij}$} The rate at which an infected agent in compartment $ I_j $ infects susceptibles in compartment $ S_i $ depends on the contact structure of a population as well as the probability that a contact between a susceptible and infectious agent leads to a transmission of the disease.
We base our contact process on data from the POLYMOD study on daily contacts in several European countries \cite{MossHens08}.
From this data we calculate a contact matrix $ C =  (c_{ij}) $ whose $ (i,j) $-th entry is the mean amount of contacts an individual in age group $ i $ has with age group $ j $; here we only consider those contacts labeled as \textit{physical}, since those are more likely to lead to viral transmission.

Let us denote by $ \beta^0_{ij} $ the rate at which a single infectious agent from age group $ j $ infects susceptible agents from age group $ i $ if no countermeasures, such as social distancing, are in place.
We model $ \beta^{0}_{ij} $ to be proportional to $ c_{ij} $ and let $ \alpha $ be the corresponding proportionality constant.
If a single infectious agent is introduced without interventions such as test, quarantine and social distancing measures in place into the otherwise completely susceptible, i.e. virgin, population, the mean amount of secondary cases he causes is the basic reproduction number 
\begin{align}
    \label{eq:r0}
    \mathcal{R}_0 \approx \sum\limits_{i,j = 1}^{n_g} N_i N_j \beta^{0}_{ij} \frac 1 {\eta} = \sum\limits_{i,j = 1}^{n_g} N_i N_j c_{ij} \frac \alpha {\eta}.
\end{align}

There is a wide variety of estimates for $ \mathcal{R}_0 $ in the literature~\cite{park_systematic_2020}, with most estimates in the interval $ [2, 3.5] $. 
We choose a value of $ \mathcal{R}_0 = 2.5 $ as early, higher estimates might be biased upwards due to imported and undetected cases. 
Fixing $ \eta^{-1}  = 6$ (see the discussion on $ \eta_A $ below) we calculate $ \alpha = 5.79 \% $ and in turn $\beta^{0}_{ij}$ from \eqref{eq:r0}:
\begin{align}
    \label{eq:beta_nominal}
    (\beta^{0}_{ij})_{1 \leq i, j \leq n_g} = \begin{pmatrix}
        0.46 & 0.48 & 0.12 \\
        0.48 & 0.63 & 0.29 \\
        0.12 & 0.29 & 0.18 
    \end{pmatrix}.
\end{align}

\paragraph{\boldmath $\gamma$} 
The rate at which latent cases become infectious is the inverse of the mean incubation time. 
This parameter is modeled age-independent and chosen to be $ 0.19 $, which corresponds to a mean incubation time of $ 5.2$ days \cite{lauer_incubation_2020}.

\paragraph{\boldmath $\pi_i^S, \pi_i^M, \pi_i^A$} These parameters denote the proportion of individuals in age group $ i $ that have severe, mild or asymptomatic course of disease. 
For Germany, the Robert Koch Institut (RKI) has published data on severity of disease progression for $ 12,178 $ cases by age-groups~\cite{schilling_vorlaufige_2020}.  
For our purposes we define a severe case to be a case that will eventually be admitted to intensive care, a mild case being one developing influenza-like symptoms, pneumonia or being admitted to hospital for other reasons. 
All other cases we classified as asymptomatic. 

Thus we obtain $ \pi_i = (\pi^S_i, \pi^M_i, \pi^A_i)$, the proportion of severe, mild and asymptomatic cases in age group $ i $, respectively, as
\begin{align*}
    \left( \pi_1, \pi_2, \pi_3 \right) = 
    \frac 1 {100}
    \begin{pmatrix}
        0.53 & 0.31 & 3.02 \\
        12.11 & 22.01 & 25.12 \\
        87.37 & 77.68 & 71.86 
    \end{pmatrix}.
\end{align*}
Observe that the oldest age group is at highest risk with $ 3.02 \% $ of infected individuals admitted to ICU. 
Also the proportion of severe cases in the youngest age group is higher than in the middle age group.
This might be explained by the fact that cases in the youngest age group are detected less frequently due to them being tested less, leading to overreporting of severe cases.

\paragraph{\boldmath $\eta^S, \eta^M, \eta^A$} 
These are the rates at which infectious individuals are removed from the infection process, if no mass-testing is implemented, i.e. if $ \theta_i = 0 $. 
For individuals with severe or mild course of disease when they develop symptoms leading to self-isolation, quarantine prescribed by a physician, or to direct hospitalization. 

One characteristic of COVID-19 is that even presymptomatic cases transmit SARS-CoV2 \cite{he_temporal_2020}.
We assume the time from being infectious to symptom onset to be $ 2 $ days after which we add $ 2 $ more days which it takes before the infectee visits a physician. Thus we choose $ \eta^S = 0.25$.

For mild progressions we assume the same mean duration from being infectious to symptomatic, though in this case individuals self-isolate, visit a physician or receive a positive test result after a mean waiting time of again $ 2 $ more days; consequently, we also set $ \eta^M = 0.25$.

For asymptomatic cases in $ I_i^A $ the only way to be removed from the infection process is by recovery from the infection. 
In \cite{woelfel_clinical_2020} positive virus samples were found in patient's throats for up to $ 8 $ days after symptom onset. 
Assuming a lower viral load for asymptomatic cases with only $ 4 $ days of potential infectiousness and adding the $ 2 $ days of presymptomatic transmission we chose $ \eta^A = 0.17 $, corresponding to a mean time of $ 6 $ days to recovery for asymptomatic cases.

\paragraph{\boldmath $\tau^O, \tau^S$} 
As we assume the testing related to the controls $ \theta_i $ to be of a random nature, tested individuals are not yet removed from the infection process. 
Instead we assume positive test results to become available after a mean delay of $ 2 $ days. 
However, severe cases may visit a physician and thus go into immediate quarantine before receiving their test result.
The latter transition occurs with rate~$ \eta^S $, and hence the faster transition occurs with rate $ \tau^S = \frac 1 2 + \eta^S = 0.75 $.

Non-severe cases that are tested, $T_i^O $, are removed if they recover naturally (with rate $ \eta^A $), or receive a positive test result, or visit a physician. 
That leads to $ \tau^O = \eta^M + \eta^A  + \frac 1 2 = 0.92$ for each age group.
\paragraph{\boldmath $\rho$} 
This parameter is the rate at which severe cases move from being in the pre-ICU state to the intensive care unit. 
This includes time spent in quarantine at home as well as time spent in the hospital in normal care while being isolated. 

In \cite{dreher_charakteristik_2020} the median time from symptom onset to being in intensive care for 50 patients was $ 9 $ days. 
As the median of an exponential distribution $ \operatorname{Exp}(\zeta) $ is $ \frac{\log 2}{\zeta} $ we choose a mean stay of $ \frac 9 {\log 2} - 2\approx 10.98 = \rho^{-1}$ days, accounting for the mean two days from symptom onset to the transition into the pre-ICU compartment.
\paragraph{\boldmath $\sigma$} 
This is the mean time spent on intensive care until discharge or death. 
According to \cite{dreher_charakteristik_2020} patients with acute respiratory distress symptoms (ARDS) spent a median amount of $ 13 $ days in intensive care and patients without ARDS spent a median amount of $ 2 $ days in intensive care. 
Of the $ 50 $ patients considered in this study, $ 24 $ were afflicted with ARDS.
Converting again between median and mean for the assumed exponential distribution yields a mean time of $\sigma^{-1} = \frac {13} {\log 2} \frac {24} {50} + \frac 2 {\log 2} \frac{ 26 } {50} = 10.5$ days spent in intensive care.
\paragraph{\boldmath $ H^{\mathrm{ICU}}_{\mathrm{max}} $}
The DIVI-Intensivregister offers daily information on the amount of free intensive care hospital beds in Germany. 
On 20 October 2020 they reported a capacity of $ 8,872 $ free beds with $ 879 $ actively treated COVID-19 patients \cite{divi_october_2020}.
We therefore round the maximal ICU-beds available for COVID-19 patients to $ 10,000 $.

\paragraph{\boldmath $ T^{\mathrm{max}}$}
In late August until the beginning of October the RKI conducted between 1 and 1.2 million weekly SARS-CoV-2 tests in Germany. 
This motivates our upper bound $ T^{\mathrm{max}} = \frac{1,200,000}{7} $ of daily tests.

\paragraph{\boldmath $ x^0 $}
We initialize our model at time $ t = 0 $ with entries of $ x^0 $ set to $ 0 $ except for those related to the susceptible, latent and infectious compartments. 
Our choice of initial values is informed by the number of active cases reported by the RKI in late march assuming the proportion of underreporting to be $ 50\% $.
We hence set the total number of infectious agents at $ t = 0 $  to 524 and the number of latent agents to 1672 distributed among the age-groups according to $ N_i $.
As we explain in Remark \ref{rem:turnpikes} our model is robust against misspecification of the initial values.

\begin{figure}[htbp!]
\centering
\includegraphics[width=\textwidth]{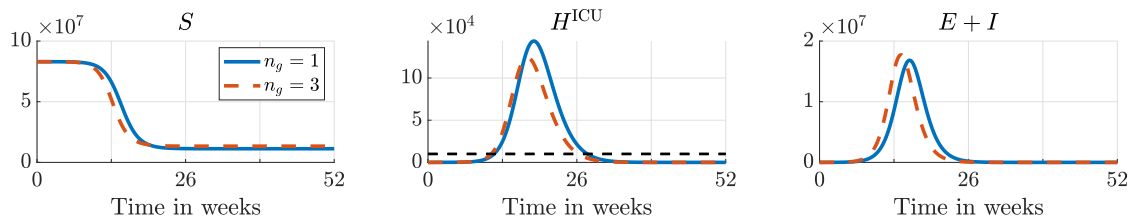}
\caption{Evolution of the pandemic without countermeasures for one and three age groups over one year; the dashed black horizontal line in the middle figure marks $ H^\mathrm{ICU}_\mathrm{max} $.}
\label{fig:sim_one_and_three_groups}
\end{figure}

Figure~\ref{fig:sim_one_and_three_groups} demonstrates the simulation capabilities of our model. Here, the course of the pandemic is visualized if no countermeasures are implemented, i.e. no social distancing ($\beta(t) = \beta^0$) and no mass-testing ($\theta(t) = 0$). %
As expected, the pandemic evolves too fast to satisfy any reasonable cap on the number of required ICUs. 
In particular, the number~$H^\mathrm{ICU}$ of required ICUs exceeds 100,000 whereas we noted above that in Germany only about 10,000 ICUs are available to treat COVID-19 patients. Therefore, countermeasures are indispensable to avoid an overload on the hospitals. %
Note that if we distinguish different age groups the pandemic evolves faster, but less ICUs are required, as the pandemic spreads mostly in the less vulnerable, younger age groups. 
Similar observations, viz. herd immunity being achieved faster in heterogeneous populations in comparison to homogeneous ones, have already been made by~\cite{britton_mathematical_2020}.

\section{Optimal testing and social distancing}\label{sec:optimal_control}
In this section, we provide information on how to keep the epidemic manageable. To this end, we formulate suitable optimal control problems (OCPs) and solve them numerically. %
Since we neither take vaccines nor re-infections into account, we consider the epidemic to be over once herd immunity is achieved, i.e. a state where the introduction of new infectious agents does not lead to an outbreak. %
Therefore, our main goal is to reach herd immunity with as few social distancing as possible while maintaining strict limits on the ICU occupancy to avoid a breakdown of the health care system. %
We call a control $u = (\beta, \theta)$ of the system~\eqref{eq:xdot=f(x,u,p)} \emph{feasible}, if $\beta_{ij}(t) \in [0, \beta_{ij}^0]$, $\theta_i(t) \geq 0$, $i,j \in \{1, \ldots, n_\mathrm{g}\}$, and
\begin{align}
	n_\mathrm{pop} \cdot \sum_{i=1}^{n_\mathrm{g}} H_i^\mathrm{ICU}(t) \; \leq \; H_\mathrm{max}^\mathrm{ICU} \notag
\end{align}
is satisfied for all~$t$.% 

A natural stopping point for simulations is when the share of susceptibles has decreased enough to ensure herd immunity even when all countermeasures are lifted completely. %
The time-dependent effective reproduction number~$\mathcal R(t)$, that is the mean number of secondary cases a primary cases will cause at time $ t $, can be used to determine whether herd-immunity has been reached: this will be the case if~$\mathcal R(t)$ is less than~$1$. %
If there is only one infected compartment, as in a simple SIR model, the latter condition is equivalent to $ \dot I (t) < 0 $. %
If there is more than one infected compartment, as in model \eqref{eq:model}, \cite{DiekHees90,DrieWatm02} have suggested to compute $ \mathcal R^{\mathrm{NGM}}(t) $, based on the so-called next-generation matrix, as a proxy for the effective reproduction number. %
Then, $ \mathcal R^{\mathrm{NGM}} (t) $ exhibits the same threshold property as~$ \mathcal R(t) $, that is $ \mathcal R^{\mathrm{NGM}} (t) < 1 $ implies herd immunity. %
Thus we use $  \mathcal R^{\mathrm{NGM}} (t) $ to check whether our simulations have reached herd immunity. %
A time horizon of two years (104 weeks) sufficed for all our simulations. %

This section is structured as follows. %
First, we verify the existence of a feasible testing strategy, i.e. without enforcing social distancing. Note that due to delays in testing, the existence of a solution is not trivial and depends on the initial value. Next, we establish an upper bound on the maximal number of tests per day and investigate to what extent social distancing is required in order to ensure feasibility. Throughout our simulations, we assume the length of one control interval to be one week. This  reflects the practical constraint that the government cannot change policies arbitrarily often but more realistically on a weekly basis. %
Throughout our simulations we use the \texttt{Matlab}-inherent sequential quadratic programming (SQP) tool to solve the OCPs. %

\subsection{Optimal testing strategy}\label{sec:optimal_testing}
Here, our goal is to maintain a hard cap on the number of required ICUs with as few tests as possible without enforcing social distancing, i.e. $\beta \equiv \beta^0$. %
To this end, we solve the OCP %
\begin{subequations}\label{theta_problem}
\begin{align}
	\min_{\theta} \quad & J_1(\theta) = \int_{t_0}^{t_f}{ T^\mathrm{tot}(t) \, \mathrm{d} t} \\
	\mathrm{s.t.} \quad & n_\mathrm{pop} \cdot \sum_{i=1}^{n_\mathrm{g}} H_i^\mathrm{ICU}(t) \leq H_\mathrm{max}^\mathrm{ICU} \quad \forall \, t \in [t_0, t_f] \label{eq:theta_problem_ICU_cap} \\
	& \dot{x}(t) = f(x(t),u(t),p), \quad \forall \, t \in [t_0, t_f], \quad x(0) = x^0 \label{eq:theta_problem_dynamics} \\
	& \beta(t) = \beta^0 \quad \forall \, t \in [t_0, t_f) \\
	& \theta_i(t) \geq 0 \quad \forall \, t \in [t_0, t_f), \, i \in \{1,2,\ldots,n_\mathrm{g}\}. \label{eq:theta_problem_box_constraints} 
\end{align}
\end{subequations}
The objective function penalizes the total number of tests over the entire time horizon~$[t_0, t_f]$ with $t_f > t_0 \geq 0$. The equality constraint~\eqref{eq:theta_problem_dynamics} captures the system dynamics while the one-sided box constraints~\eqref{eq:theta_problem_box_constraints} ensure that the testing rates cannot be negative. 

\begin{figure}[htbp!]
\centering
	\includegraphics[width=\textwidth]{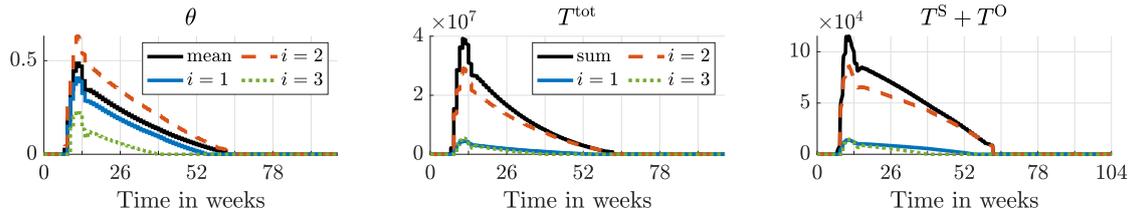} 
	\caption{Optimal testing strategy for three age groups over two years.}% 
	\label{fig:results_theta_controls}
\end{figure}

\begin{figure}[htbp!]
\centering
	\includegraphics[width=\textwidth]{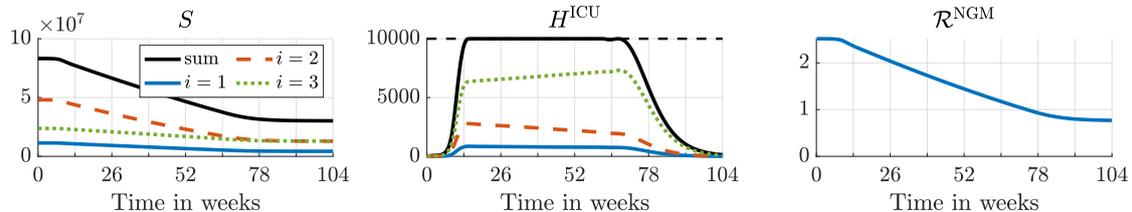}
	\caption{Evolution of the epidemic corresponding to the optimal testing strategy visualized in Figure~\ref{fig:results_theta_controls} and reproduction number based on next-generation matrix. The dashed black line in the middle plot depicts the maximal capacity of ICUs.}
	\label{fig:results_theta_states}
\end{figure}

Figure~\ref{fig:results_theta_controls} depicts the optimal controls as well as the total number of tests and the number of detected cases per day while Figure~\ref{fig:results_theta_states} shows the impact on the evolution of the epidemic. %
Here, we computed the effective reproduction number~$\mathcal{R}^\mathrm{NGM}(t)$ at each time step, demonstrating that we reached herd immunity.

We observe that there exists a testing strategy that ensures feasibility, which was not obvious from the outset because of the assumed delays. %
In particular, the bound~\eqref{eq:theta_problem_ICU_cap} is active once it's reached, i.e. $H^\mathrm{ICU} \equiv H_\mathrm{max}^\mathrm{ICU}$, and becomes inactive when the number of susceptible people falls below a certain threshold and $\mathcal{R}^\mathrm{NGM}(t) < 1$ indicating the onset of herd immunity. %
\begin{remark}\label{rem:turnpikes}
	The steady state in~$H^\mathrm{ICU}$ suggests that problem~\eqref{theta_problem} satisfies the so-called turnpike property~\cite{Grue16}. Typically, turnpikes indicate the optimal operating state of a system. These are steady states at which the running costs are minimized. In our example, since we do not penalize the number of required ICUs, the best strategy is to stay at the upper bound while saving tests. Once the objective function value is zero the system leaves the state eventually. In particular, regardless of the initial value, the system is steered towards this optimal operating point. As a consequence, a rough estimation suffices as initial guess for our simulations. %	
A rigorous analysis of these turnpikes, however, is left for future research. 
\end{remark}

However, these results are only of theoretical interest, since this optimal testing strategy would be prohibitively expensive and might not even be implementable at all. For instance, regarding Figure~\ref{fig:results_theta_controls}, one observes that the mean testing rate reaches about~$0.5$, which corresponds to being tested every two days on average. Moreover, the total number of tests per day %
required for this approach is more than 12,000,000 %
on average (over 65~weeks), %
compared to $ T^{\mathrm{max}} \approx 170,000$ daily tests which are currently conducted in Germany. Note that, even with this enormous testing effort, the number of detected cases, $T^\mathrm{S} + T^\mathrm{O}$, is rather small since the number of infectious individuals is small compared to the total population. %

In conclusion, mass-testing alone currently does not suffice to maintain hospitalization caps in reality. These arguments support the government's decision to introduce additional measures like social distancing and hygiene concepts. 
However, cheap rapid test kits might change the situation favorably as they could be made widely available, self-administered while giving immediate test results.

In the following subsections, we enforce $ T^\mathrm{max} $ as an upper bound on the amount of daily tests.
Under this additional constraint we then determine the minimal amount of social distancing required to reach herd immunity.
The success of such measures depends on the acceptance and thus compliance by the general population.

\subsection{Optimal homogeneous social distancing}
In a first step, we determine an optimal social distancing strategy by penalizing the deviation of~$\beta$ from~$\beta^0$ equally over all age groups. %
This might increase acceptance in the general population due to the (perceived) fairness of such measures: everyone is treated equally and contacts are reduced by the same proportion for everyone. 
In reality such strategies may be hard to conceive as different measures affect the age groups differently, i.e. closing schools and nurserys affects those in the lowest age group the most while leaving the oldest age group unaffected.
Nevertheless a mixture of many different non-pharmaceutical measures may be able to achieve such a reduction in contacts.

We introduce a time-varying factor $\delta = \delta(t)$ describing the amount of social distancing that is implemented. Moreover, we choose to penalize the $\ell^2$ deviation of this control input from $\delta = 1$ in the objective function in order to smooth the optimal control. For instance, penalizing the $\ell^1$ deviation yields bang-bang controls, i.e. the optimal solution jumps back and forth between the two extremal options: no contact restrictions and lock down (simulations not shown). %
Therefore, we %
determine an optimal homogeneous social-distancing policy by solving% 
\begin{subequations}\label{delta_problem}
\begin{align}
	\min_{\theta,\delta} \quad & J_2(\theta, \delta) = \int_{t_0}^{t_f}{ \big(1 - \delta(t)\big)^2 + \kappa \sum_{i=1}^{n_\mathrm{g}} \theta_i(t) \, \mathrm{d} t} \\
	\mathrm{s.t.} \quad & n_\mathrm{pop} \cdot \sum_{i=1}^{n_\mathrm{g}} H_i^\mathrm{ICU}(t) \leq H^\mathrm{ICU}_\mathrm{max} \quad \forall \, t \in [t_0, t_f] \\
	& \beta_{ij}(t) = \delta(t) \beta_{ij}^0 \quad \forall \, i,j \in \{1, \ldots, n_\mathrm{g}\} \\
	& \dot{x}(t) = f(x(t),u(t),p) \quad \forall \, t \in [t_0, t_f], \quad x(0) = x^0 \label{eq:delta_ode} \\
	& T^\mathrm{tot}(t) \leq T^\mathrm{max} \quad \forall \, t \in [t_0, t_f], \label{eq:delta_theta_max} \\
	& \delta(t) \in [0, 1] \quad \forall \, t \in [t_0, t_f) \\ 
	& \theta_i(t) \geq 0 \quad \forall \, t \in [t_0, t_f), \, i \in \{1,2,\ldots,n_\mathrm{g}\}. \label{eq:delta_theta_min} 
\end{align}
\end{subequations}
Note that we allow to distribute the tests among the age groups by not fixing~$\theta_i$, but enforcing~\eqref{eq:delta_theta_max} and~\eqref{eq:delta_theta_min}. Furthermore, we introduce a regularization term with weight $\kappa = 10^{-5}$. %
The choice of~$\kappa$ is based on simulations. 
In contrast to~\eqref{theta_problem}, we always find a feasible solution of~\eqref{delta_problem} if the epidemic has not yet evolved too far. More precisely, by choosing $\delta = 0$, which corresponds to a complete lockdown, we are (theoretically) able to stop the spread. Therefore, if the initial number of people with a severe course of infection is sufficiently low, the upper bound on the number of ICUs will not be violated. 

\begin{figure}[htbp!]
\centering
\includegraphics[width=\textwidth]{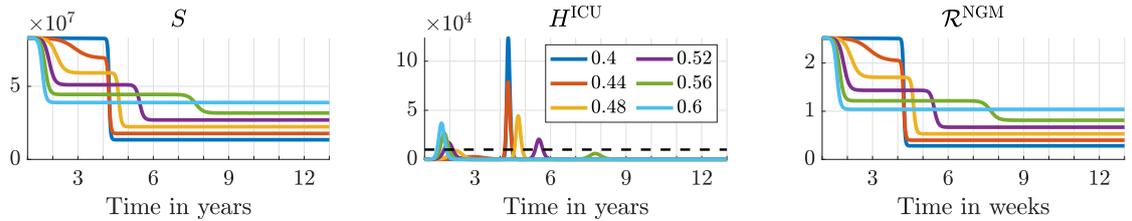}
\caption{Solution for fixed $\delta$ for three years and complete lift of restrictions afterwards.}
\label{fig:constant_delta_let_go}
\end{figure}

A highly fluctuating social distancing strategy may lead to low acceptance in the general population, because people have to adapt to new rules every few weeks.
Thus before we solve~\eqref{delta_problem} let us have a look at what happens if we consider a constant value for~$\delta$ over time, i.e. a social distancing strategy without fluctuations. %
\begin{figure}[htbp!]
\centering
	\includegraphics[width=\textwidth]{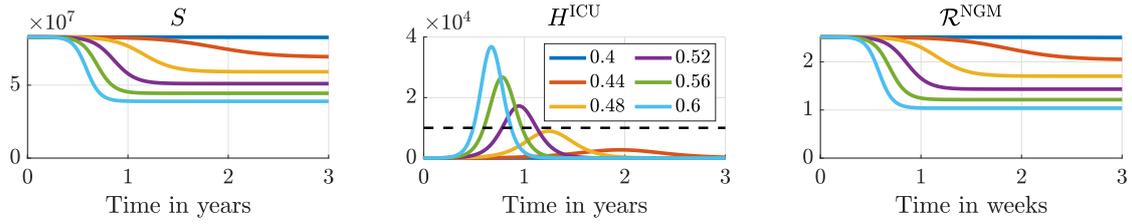}
	\caption{The impact of constant contact reduction rate~$\delta \in [0.4, 0.6]$ on the speed of evolution of the epidemic and on the number of required ICUs.}
	\label{fig:constant_delta}
\end{figure}
Figure~\ref{fig:constant_delta} (left) shows that fewer contacts result in a longer time for the epidemic to abate on the one hand, but a lower number of total infections within the considered time horizon on the other hand. Moreover, Figure~\ref{fig:constant_delta} (middle) visualizes that quite strict social distancing is needed in order to meet the ICU capacities. The maximal value of~$\delta$ to stay feasible is~$0.487$, i.e. contacts needed to be more than halved over three years. Furthermore, once we lift the restrictions, see Figure~\ref{fig:constant_delta_let_go}, there might be another outbreak. %
In particular, the stronger the restrictions were in the beginning, the stronger the second outbreak will be. 
Therefore, it is essential to establish herd immunity before lifting all restrictions, and to adapt the policy over time. %

\begin{figure}[htbp!]
\centering
	\includegraphics[width=\textwidth]{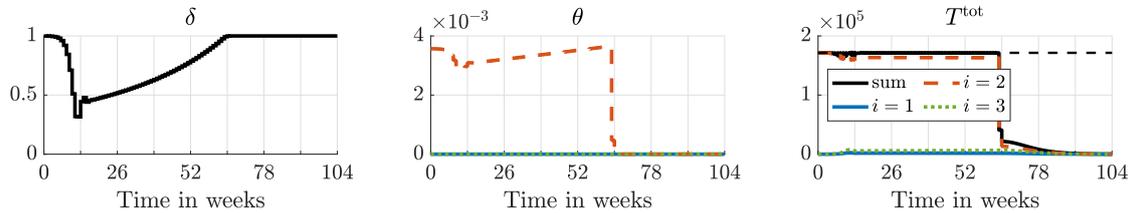} 
	\caption{Optimal combination of testing and (homogeneous, time-varying) social distancing for three age groups over two years. The dashed black line in the second subplot depicts the upper bound on the total number of tests per day.} %
	\label{fig:results_delta_controls}
\end{figure}

\begin{figure}[htbp!]
\centering
	\includegraphics[width=\textwidth]{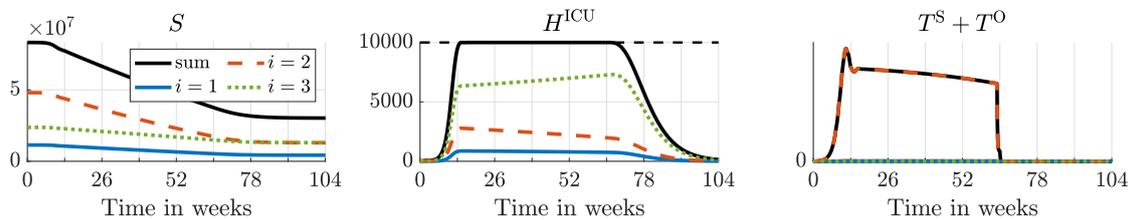}
	\caption{Evolution of the epidemic based on the controls visualized in Figure~\ref{fig:results_delta_controls}. The dashed black line in the plot of $H^\mathrm{ICU}$ depicts the upper bound on the number of available ICUs.}
	\label{fig:results_delta_states}
\end{figure}

A visualization of the optimal solution of~\eqref{delta_problem} can be found in Figures~\ref{fig:results_delta_controls} and~\ref{fig:results_delta_states}. %
As mentioned above, the bound on~$H^\mathrm{ICU}$ is not violated. Since the weight~$\kappa$ is chosen sufficiently small, the upper bound on the total number of tests per day is active as long as the upper bound on~$\delta$ is not. However, note that not all age groups are tested equally. More precisely, only the middle-aged group is tested at all. The reason is that this group is the largest ($N_2 > N_1 + N_3$) and has the highest contact rates (c.f.~\eqref{eq:beta_nominal}) and therefore, contributes more to the spread of the epidemic than the other groups. Furthermore, we observe that the social distancing policy has to be quite strict in the beginning. In particular, $\min_t \delta(t) \approx 0.3$ which corresponds to a reduction of average contacts per person by $70\%$. 
However, this can be qualitatively compared to the measures taken in Germany 
starting in mid March 2020 when contacts were reduced by school and restaurant closures as well as other contact restrictions. %

In conclusion, social distancing is an effective tool to keep the epidemic manageable. Comparing the results of~\eqref{delta_problem} to the simulations with constant~$\delta$ we see that a (partial) lockdown appears inevitable. However, our simulations suggest to let the epidemic evolve for a few weeks, then enforce a contact reduction down to approximately 30\% for 2-4 weeks before slowly lifting the restrictions over the next 12~months until herd immunity is achieved. 

\subsection{Age-dependent social distancing}
The constraint that contacts are reduced by the same proportion for each age group is restrictive and it is plausible that more efficient solutions exist when contact reductions are distributed differently across age groups. 
One reason to consider such a strategy is that it may be more efficient at stopping the spread of the epidemic; as mentioned above the middle-age group is the driver of the epidemic while the oldest age group consists of the most vulnerable individuals.
In any case, such an age-differentiated social distancing strategy needs to be accepted by the whole population to be successful.

Hence, we improve the social distancing policy computed above by allowing it to depend on age. Given the solution~$(\theta^\star,\delta^\star)$ of~\eqref{delta_problem}, %
we solve the OCP %
\begin{subequations}\label{beta_problem_b}
\begin{align}
	\min_{\theta, \beta} \quad & J_3(\theta, \beta) = \int_{t_0}^{t_f}{ \sum_{i,j=1}^{n_\mathrm{g}} N_i N_j \big(\beta_{ij}(t) - \beta_{ij}^0\big)^2 + \kappa \theta_i(t) \, \mathrm{d} t} \\
	\mathrm{s.t.} \quad & n_\mathrm{pop} \cdot \sum_{i=1}^{n_\mathrm{g}} H_i^\mathrm{ICU}(t) \leq H^\mathrm{ICU}_\mathrm{max} \quad \forall \, t \in [t_0, t_f] \\ 
	& \dot{x}(t) = f\big(x(t), u(t), p\big) \quad \forall \, t \in [t_0, t_f], \quad x(0) = x^0 \\
	& T^\mathrm{tot}(t) \leq T^\mathrm{max} \quad \forall \, t \in [t_0, t_f] \\
	& \theta_i(t) \geq 0 \quad \forall \, t \in [t_0, t_f), \, i \in \{1,2,\ldots,n_\mathrm{g}\} \\
	& \beta_{ij}(t) \in [\beta_{ij}^\mathrm{min},\beta_{ij}^0] \quad \forall \, t \in [t_0, t_f), \quad i,j \in \{1,2,\ldots,n_\mathrm{g}\}. 
\end{align}
\end{subequations}
Here, we use $\delta^\star$ %
to define $\beta_{ij}^\mathrm{min} = \min_t \delta^\star(t) \beta_{ij}^\mathrm{nom}$, %
i.e. the lower bound on~$\beta$ in~\eqref{beta_problem_b} is the worst case %
of~\eqref{delta_problem}. %
Therefore, no one is treated worse than when applying homogeneous social distancing. %
Note that $(\theta^\star, \beta)$ with $\beta = \delta^\star \beta^0$ is feasible for problem~\eqref{beta_problem_b}. %
As in~\eqref{delta_problem} we penalize testing as soon as $\beta(t) = \beta^0$ holds. %

\begin{figure}[htbp!]
\centering
	\includegraphics[width=\textwidth]{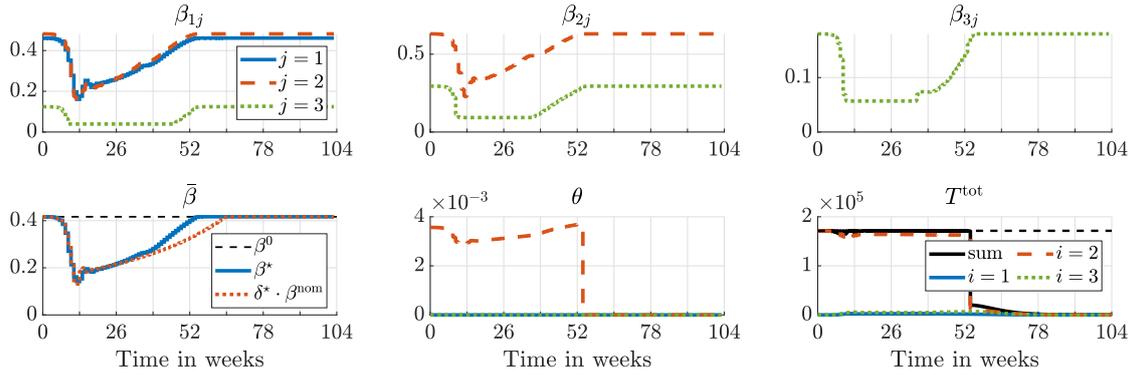} 
	\caption{Optimal age-dependent social-distancing strategy for three age groups over two years.
	}
	\label{fig:results_beta_controls_b}
\end{figure}

\begin{figure}
\centering
	\includegraphics[width=.67\textwidth]{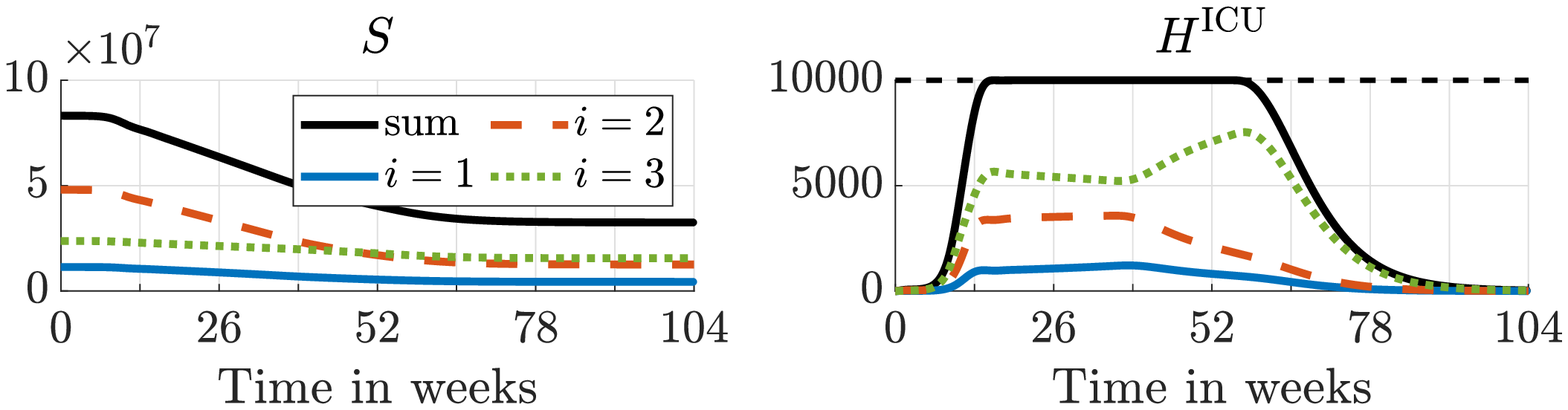}
	\caption{Evolution of the compartments associated with controls depicted in Figure~\ref{fig:results_beta_controls_b}.}
	\label{fig:results_beta_states_b}
\end{figure}

Results for~\eqref{beta_problem_b} can be found in Figures~\ref{fig:results_beta_controls_b} and~\ref{fig:results_beta_states_b}, where
\begin{align}
	\bar{\beta}(t) \; = \; \sum_{i,j=1}^{n_\mathrm{g}} N_i N_j \beta_{ij}(t) \notag
\end{align} 
describes the average number of contacts per person and day in a heterogeneous population. %
Here, we used~$\kappa = 10^{-5}$. %
The corresponding value for $ \beta^0 $ is $ \bar\beta^0 =  0.4167$. %
This allows to compare the solution $ \beta_{ij}(t) $ with $ \bar\beta^0 \delta^\star(t) $ obtained from \eqref{delta_problem}. 
Similar to the solution of~\eqref{delta_problem}, the upper bound on testing is active most of the time, while essentially only the middle-aged group is tested. %
The social distancing measures %
are less restrictive than for~\eqref{delta_problem} which makes compliance with the measures more likely. %
However, the measures could be perceived as unfair, since the contacts of the oldest age group are restricted most. %
Moreover, the contacts of the middle-age group are least restricted. Therefore, the working class would be allowed to go to work which benefits the economy.

In conclusion, social distancing is crucial to avoid an overload on the hospitals. In addition, testing middle-aged people helps to reduce the required amount of social distancing. Furthermore, all presented strategies support a lock down a few weeks into the epidemic, which is followed by lifting the restrictions step by step until herd immunity sets in. %
Age-differentiated social distancing might be hard to argue for, but it helps to end the epidemic several months earlier and, therefore, support the economy. 

\subsection{Short-term decision making}
The control strategies derived in the previous subsections provide rough guidelines for how the epidemic can be controlled. However, from a decision maker's perspective, it will be hard to argue for policies taking effect in the far future. In particular, there are many uncertainties that might affect the performance of the control strategy over the time span of two years, and hence the control strategy needs to be adjusted over time. %
Therefore decision need to be revised constantly adapting to the changing conditions during the epidemic's progress. %
Model predictive control (MPC) provides a state of the art methodology to tackle such problems. %
The basic idea of MPC is to consecutively solve a series of OCPs over a smaller horizon of~$K$ control intervals rather than solving a single OCP over the whole horizon. More precisely, only the first part of the optimal control derived by solving such an auxiliary OCP is implemented. Then, the time window is shifted, and the procedure is repeated based on updated measurements. 
For a detailed introduction to MPC we refer to~\cite{RawMD19}. %
Here, we tackle~\eqref{beta_problem_b} via MPC; the earlier problems can be treated analogously. The MPC scheme for~\eqref{beta_problem_b} is summarized in Algorithm~\ref{alg:mpc}. %
\begin{algorithm}%
\caption{MPC scheme for solving~\eqref{beta_problem_b}}
{\bf Input:} Prediction horizon length~$K$, length of control interval~$\Delta t$. Set time $t = t_0$. \\
{\bf Repeat:}
\begin{enumerate}
	\item Obtain current states $\hat{x} = x(t)$.
	\item Determine optimal solution $u^\star : [t,t+K\Delta t) \to \mathbb{R}^9$ of~\eqref{beta_problem_b} on $[t,t+K\Delta t)$ with $x^0 = \hat{x}$.
	\item Implement $u^\star|_{[t,t+\Delta t)}$. Increment time $t \leftarrow t + \Delta t$.
\end{enumerate}
\label{alg:mpc}
\end{algorithm}

\begin{figure}%
\centering
\includegraphics[width=\textwidth]{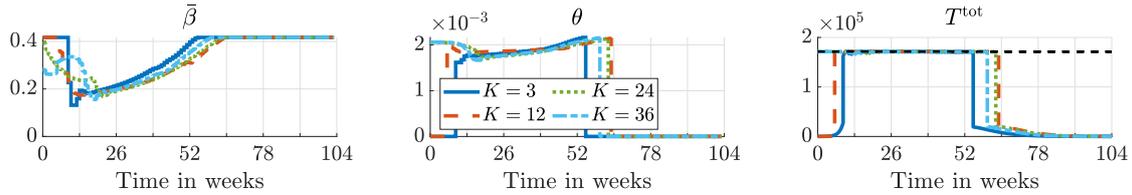}
\caption{Optimal control for solving~\eqref{beta_problem_b} in closed loop for varying prediction horizon length. 
For the sake of readability, we depicted average values of~$\theta$ and sums of~$T^\mathrm{tot}$ over the age groups.}
\label{fig:beta_controls_CL}
\end{figure}

\begin{figure}%
\centering
\includegraphics[width=.67\textwidth]{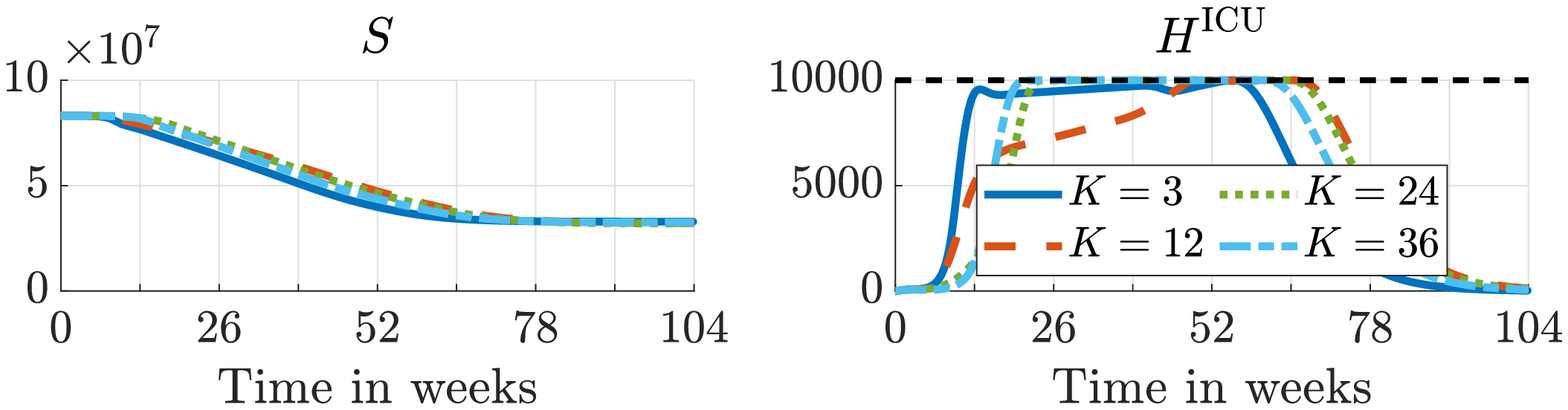}
\caption{Evolution of the epidemic based on the controls depicted in Figure~\ref{fig:beta_controls_CL}. For the sake of readability, only the sum over the age groups is visualized.}
\label{fig:beta_states_CL}
\end{figure}

Results based on varying prediction horizon lengths can be found in Figures~\ref{fig:beta_controls_CL} and~\ref{fig:beta_states_CL}. %
The basic structure of both the optimal control and the associated states is comparable to the open-loop solution presented in the previous subsection. Therefore, we stopped the simulations after one and a half years. The length of the prediction horizon affects mainly the optimal social distancing policy. In particular, the larger the prediction horizon, the less social distancing is needed in total. More precisely, for bigger~$K$, we implement a slightly stricter lockdown but can start it later and relax it earlier. Furthermore, the larger~$K$ gets the closer the optimal solution is to the open-loop solution.
In particular, the MPC solutions qualitatively resemble the open-loop solution: after an early lockdown, social distancing is slowly lifted.

For $K = 3$, the ICU capacity reaches its upper limit earlier due to the laissez-faire policy in the beginning. However, this constraint also becomes inactive earlier.
For even shorter prediction horizons recursive feasibility cannot be guaranteed, i.e. the ICU cap might be violated (simulations not shown).

\subsection{Impact of upper bounds on number of tests and ICUs}
So far, we assumed both the upper bounds on the number of tests per day, and on the number of ICUs to be fixed at our chosen values.
In practice, these conditions may change: free ICU capacity might exhibit seasonal patterns and the number of possible tests per day depends on infrastructure and available personnel. 
In addition, varying the upper bounds is useful to illustrate the benefits of increased testing and higher ICU capacities.
In this subsection, we investigate the impact of these parameters on the optimal social distancing policy numerically.

\begin{figure}[htbp!]
\centering
\includegraphics[width=\textwidth]{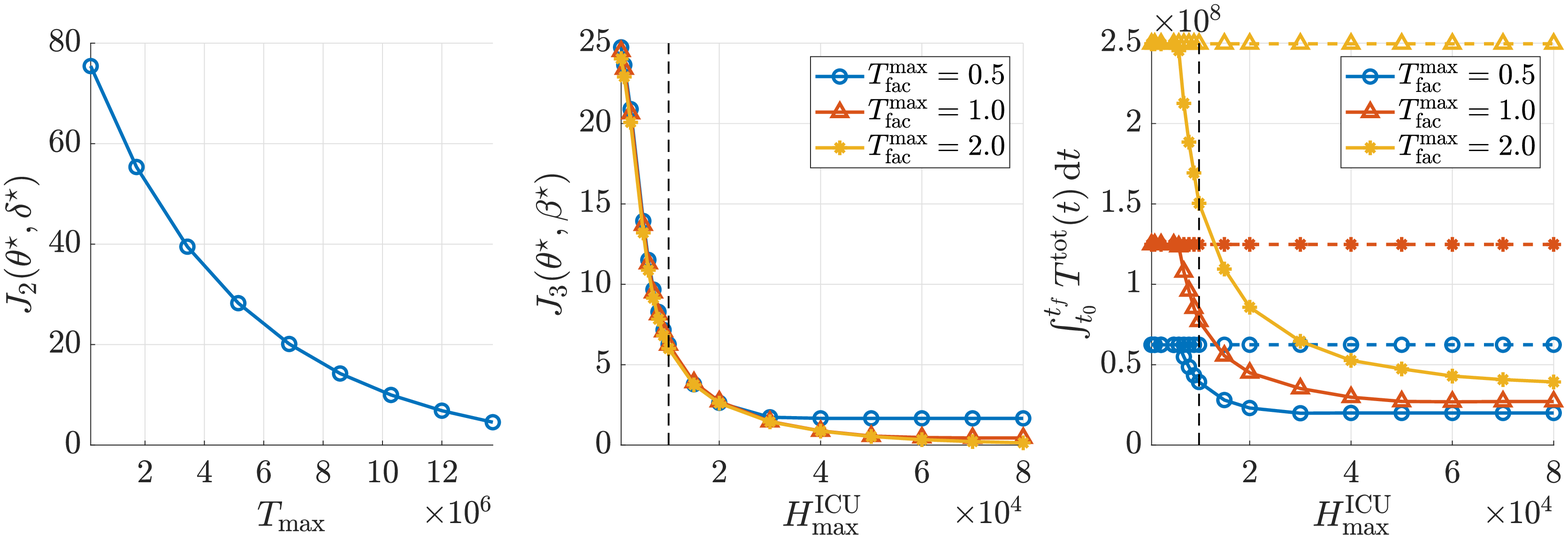}
\caption{Impact of $T^\mathrm{max}$ on social-distancing costs (left) and of $H^\mathrm{max}$ on both social-distancing costs (middle) and testing (right). In the last two subfigures the currently available number of ICUs in Germany is highlighted by a vertical dashed line. The dashed horizontal lines in the right-most figure indicate the total testing capacities over the entire simulation horizon. Factor of modification of~$T^\mathrm{max}$ denoted by~$T^{\mathrm{max}}_{\mathrm{fac}}$.}
\label{fig:pareto_all}
\end{figure}

First, we study the effect~$T^\mathrm{max}$ has on the social distancing by solving~\eqref{delta_problem} via MPC, see Figure~\ref{fig:pareto_all} (left). %
As pointed out in the previous subsection, the prediction horizon length affects the start and end time of measures as well as its peak %
(simulations not shown).
In addition, increasing~$T^\mathrm{max}$ by some factor $T^{\mathrm{max}}_{\mathrm{fac}} \geq 1$ shifts the whole $\delta$ curve upwards, i.e., as expected, the more tests are available, the less social distancing is required. Furthermore, Figure~\ref{fig:pareto_all} (left) visualizes the impact of~$T^\mathrm{max}$ on the objective function value of~\eqref{delta_problem}. %

\begin{figure}[htbp!]
\centering
\includegraphics[width=.67\textwidth]{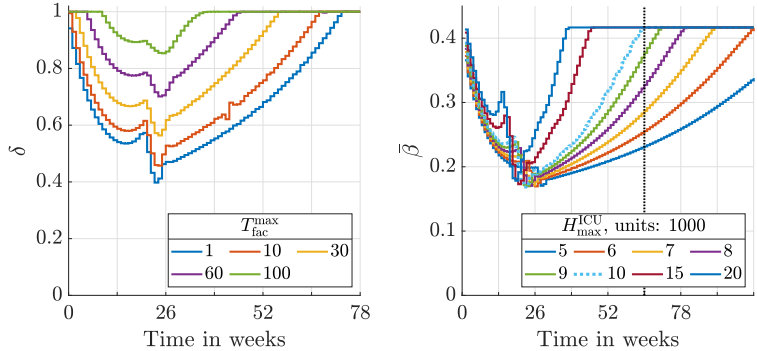}
\caption{Impact of the available number~$H^\mathrm{ICU}_\mathrm{max}$ of ICUs and the prediction horizon~$K$ on the average social distancing. The dotted cyan line refers to the number of currently available ICUs in Germany. The vertical dotted black line marks the end of social distancing measures for that setting.}
\label{fig:beta_over_time_Hmax}
\end{figure}

Second, we investigate the impact of the number of ICUs on the optimal solution of~\eqref{beta_problem_b}. Results can be found in Figures~\ref{fig:pareto_all} and~\ref{fig:beta_over_time_Hmax}. For the simulations in Figure~\ref{fig:pareto_all} (middle and right) we used MPC with prediction horizon $K=12$ weeks. Figure~\ref{fig:pareto_all} (middle) clearly shows that the number of available ICUs directly affects the cost function value. While for a small value of~$H_\mathrm{max}^\mathrm{ICU}$, every additional ICU contributes, for large values, a saturation seems to take place. In particular, doubling the current number of available ICUs does help, but the benefit becomes negligible when increasing it further. These phenomena are almost unaffected by doubling or halving~$T^\mathrm{max}$. However, when there are not enough ICUs, then the upper bound on $T^\mathrm{tot}$ is always active, see Figure~\ref{fig:pareto_all} (right), where~$T^\mathrm{tot}$ is at its maximum value all the time. %
Moreover, an increase in the number of ICUs clearly leads to a reduction in the social distancing measures, as can be seen in Figure \ref{fig:beta_over_time_Hmax}.

In summary, increasing test capacities and/or ICU capacities helps to reduce measures like social distancing. However, the impact of the number of available ICUs appears to be much stronger. 
Nonetheless the qualitative shape of the solutions over time is not affected by varying these constraints.

\section{Conclusions and outlook}\label{sec:conclusions}
In this paper, we demonstrated how mitigation of the COVID-19 epidemic can be achieved by a combination of age-stratified testing and social distancing measures while avoiding a breakdown of the health care system.
We showed that in our compartmental model mass testing alone is insufficient to achieve this goal, as it would require unrealistic testing capacities.

As a remedy, we designed optimal social distancing strategies with a focus on applicability and acceptance in the general population, i.e. strategies with slowly changing contact reductions.
The resulting social distancing measures imitate the measures actually taken in Germany, but are lifted at a much slower pace.
Age-differentiated contact reductions may improve upon these results as they yield qualitatively similar social distancing strategies and prioritize relaxing restrictions for the work-force and children.

To model the process of policy making more realistically, we used MPC which allows to adapt to deviations from the envisioned course of the epidemic by solving the optimal control problem repeatedly.
Our analysis reveals that longer prediction horizons allow for faster lifting of restrictions although long-term predictions may be infeasible in practice.
Additionally we showed that the amount of available intensive care units is a key factor influencing the required amount of social distancing.

We believe that our model with the chosen parameters reflects reality sufficiently well to provide qualitatively valid insight on how testing and social distancing can control the spread of SARS-CoV-2. 
We learned that mass testing alone is, assuming realistic testing capacities, not sufficient to avoid a breakdown of the health care system in Germany. % 
To prevent this, one has to implement strict contact reductions early on, which, ideally, should then be eased slowly. %
If one allows these reductions to vary by age, one is able to relax restrictions for the (working) middle age group, at the cost of reducing contacts of the more vulnerable older population. %

While short-term planning of measures is unable to control the exponential growth of cases, medium-term planning produces strategies that, qualitatively, do not differ from optimal ones while being flexible enough to adapt to new circumstances. %
Finally, as expected, the number of available intensive care units dictates how fast herd-immunity can be reached and how much total social distancing is necessary. %

However, we caution the reader against interpreting these results in a quantitative way, as our model has not been devised to produce precise predictions.
Similarly, we want to stress that we do not provide concrete policies to implement, as the impact of particular countermeasures on $ \beta $ is not easily quantified.

Concerning other influences on the epidemic's evolution, note that we have not yet considered vaccinations nor re-infections, both of which could be included in our model without difficulties, if parameters are available to model them. 
As our model is based on ODEs, interaction effects such as contact tracing cannot be included. 
Agent-based (stochastic) models are able to handle these critical effects and could be seen as a natural extension of our (deterministic) compartmental model. 
To solve the resulting stochastic optimal control problems would then require more sophisticated techniques, however.

	% Acknowledgements
\section*{Acknowledgments}
We thank Kurt Chudej (University of Bayreuth) for insights on modelling pandemics and Manuel Schaller (TU Ilmenau) for fruitful discussions on optimal control and the turnpike property.

	% Bibliography
\bibliographystyle{plain}
\bibliography{references}
\end{document}